\documentclass[11pt]{article}

\usepackage{amsmath} \usepackage{amsthm}
\usepackage{amsfonts}\usepackage{amssymb}
\usepackage{graphicx}
\usepackage{color}
\newtheorem{theorem}{Theorem}

\newtheorem{lemma}{Lemma}

\newcommand{\eps}{\varepsilon}
\def\ep{\epsilon}

\newcommand{\indic}[1]{\mathbf{1}_{\{#1\}}}

\def\Ee{E}
\def\Pr{P}

\def\Zz{\mathbf{Z}}

\newcommand{\half}{\frac1{2}}

\newcommand{\keywords}{\textbf{Keywords}\ }

\def\mn{\medskip\noindent}
\def\bn{\bigskip\noindent}
\def\beq{\begin{equation}}
\def\eeq{\end{equation}}
\def\beqa{\begin{eqnarray}}
\def\eeqa{\end{eqnarray}}
\def\beqax{\begin{eqnarray*}}
\def\eeqax{\end{eqnarray*}}
\def\sqz{\kern -0.2em}

\def\adj{d_{\text{adj}}}

\begin{document}

\title{Limiting behavior of the distance of a random walk}

%\title{Is the Cayley distance smooth ?}

\author{Nathana\"el Berestycki$^{1}$ and Rick Durrett$^2$}
\date{April 5, 2006}
\maketitle

\centerline{\textbf{Abstract}}

\mn This investigation is motivated by a result we proved recently
for the random transposition random walk: the distance from the
starting point of the walk has a phase transition from a linear
regime to a sublinear regime at time $n/2$. Here, we study three
new examples. It is trivial that the distance for random walk on
the hypercube is smooth and is given by one simple formula. In the
case of random adjacent transpositions, we find that there is no
phase transition even though the distance has different scalings
in three different regimes. In the case of a random 3-regular
graph, there is a phase transition from linear growth to a
constant equal to the diameter of the graph, at time $3\log_2 n$.

\vfill \noindent \keywords{random walk, phase transition,
adjacent transpositions, random regular graphs}

\bn 1. University of British Columbia. Room 121 -- 1984,
Mathematics Road. Vancouver, BC, Canada, V6T 1Z2. Ecole Normale
Sup\'erieure, D.M.A. 45, rue d'Ulm, 75005 Paris, France.

\bn 2. Department of Mathematics, Malott Hall, Cornell University,
Ithaca, NY 14853, U.S.A. Both authors are partially supported by a
joint NSF-NIGMS grant DMS-0201037.

\clearpage

\section{Introduction}
Let $X^n_t$ be the continuous time random transposition random
walk on $n$ markers. This means that at rate 1, we change the
current permutation by performing a transposition of two randomly
chosen elements. Let $D^n_t$ be the distance of $X^n_t$ from its
starting point, i.e., the minimal number of transpositions
necessary to change $X^n_t$ into $X^n_0$. The main result of
Berestycki and Durrett \cite{bd} is that $D^n_t$ has a phase
transition at time $n/2$ as $n\to\infty$. Writing $\to_p$ for
convergence in probability.

\mn \textbf{Theorem 0.}\emph{ Let $t>0$. As $n \to\infty$
$n^{-1}D^n_{tn} \to_p f(t)$
where $f(t)$ is defined by:
\begin{equation*}
f(t)=
\begin{cases}
t &\text{ for $t\le 1/2$} \\
 1 - \sum_{k=1}^{\infty}\frac1{2t} \frac{k^{k-2}}{k!}(2te^{-2t})^k <t
 & \text{ for $t>1/2$}
\end{cases}
%\right.
\end{equation*}
}

Having seen this result, it is natural to ask in what situations
is the asymptotic behavior of the distance from the starting point
non smooth. We begin with a trivial example.

\subsection{Random walk on the hypercube}

Let $X^n_t$ be the random walk on the hypercube $\{0,1\}^n$ that jumps at rate
1, and when it jumps the value of one randomly chosen coordinate is changed. By considering
a version of the chain  that jumps at rate 2, and when it jumps the new coordinate
takes on a value chosen at random from $\{0,1\}$ it is easy to see that, when $n=1$,
$$
P_0( X^1_t = 1 ) = (1-e^{-2t})/2
$$
Let $D^n_t$ be the distance from $X^n_t$ to $X^n_0$, i.e., the number of
coordinates that disagree. Since the coordinates in the continuous time change are independent
it follows easily from this that

\begin{theorem}
As $n\to\infty$, $n^{-1} D^n_{nt} \to  (1-e^{-2t})/2$ in probability.
\end{theorem}

\subsection{Random adjacent transpositions}

Let $X^n_t$ be the continuous time random adjacent transposition
random walk on $n$ markers. Here we are thinking of $X_t(j)$ as
the location of particle $j$, but the dynamics are easier to
formulate in terms of $Y_t(i) = X_t^{-1}(i) $ which is the
particle at location $i$. At rate 1, we change the permutation by
picking $1\le i\le n-1$ at random and exchanging the values of
$Y^n_t(i)$ and $Y^n_t(i+1)$. Without loss of generality we can
suppose $X^n_0$ is the identity permutation $I$. The distance from
a permutation $\sigma$ to $I$, i.e., the minimum number of
adjacent transpositions needed to build $\sigma$, is given by the
following convenient formula
\begin{equation}
\label{inv}
\adj(\sigma)=\text{Inv}(\sigma):=\#\{1\le i<j \le n:
\sigma(i)>\sigma(j)\}
\end{equation}
$\text{Inv}(\sigma)$ is called the number of inversions of
$\sigma$. If we view the set of permutations $\mathcal{S}_n$ of
$\{1,\ldots,n\}$ as a graph where there is an edge between
$\sigma$ and $\sigma'$ if and only if $\sigma'$ can be obtained
from $\sigma$ by performing an adjacent transposition (in the
sense defined above), or vice-versa, then $X_t$ has the law of
simple random walk on this graph and $\adj(X_t)$ is the length
induced by the graph distance of the shortest path between the
current state of the walk, $X_t$, and its starting point, the
identity.

\medskip Erikkson et al. \cite{erikkson+} and later Eriksen
\cite{eriksen} considered the problem of evaluating the distance
for the discrete time chain $\hat X^n_k$. Relying heavily on
formula (\ref{inv}) they were able to carry some explicit
combinatorial analysis, to obtain various exact formulae for this
expected distance, such as this one:
\begin{equation}
\label{erikkson}
E\adj(\hat X^n_{k})= \sum_{r=0}^k\frac{(-1)^r}{n^r}\left[\binom{k}{
r+1}2^rC_r+4d_r\binom{k}{r}\right]
\end{equation}
where $C_r$ are the Catalan numbers and $d_r$ is a less famous non-negative integer
sequence, which they define explicitly.

While this formula is exact, it is far from obvious how to extract
useful asymptotics from it. We will take a probabilistic approach
based on the formula
$$
D^n_t = \adj(X^n_t)=\sum_{i<j}\indic{X^n_t(i)>X^n_t(j)}
$$
If $1\le i\le n$ is fixed, the trajectory
$X^n_t(i)$ of the $i^{\text{th}}$ particle is a continuous time simple random walk on
$\{1,\ldots,n\}$ starting at $i$ with jumps at rate $2/(n-1)$
and reflecting boundaries at 1 and $n$ that cause the particle to stay put
with probability 1/2.

Two such trajectories, say those of particles $i$ and $j$ with
$i<j$, move by the nearest neighbor stirring process on
$\{1,\ldots,n\}$ (which for indistinguishable particles produces
the simple exclusion process). When the particles are not
adjacent, they perform independent simple random walks. When they
are adjacent, the only things that can happen are an exchange of
the two particles, or one of them moves away from the other. As
the reader can probably guess, and Durrett and Neuhauser
\cite{dn94} have proved on $\mathbf{Z}$, when $n$ is large the
random walks behave as if they are independent.

At small times, the behavior of the distance is messy but smooth.
For an integer $x \ge 0$, let $T^x$ denote the hitting time of the
level $x$ by a {rate} 4 random walk on $\Zz$ starting at $0$. Let
$\{Y(t), t \ge 0\}$ and $\{Y'(t), t \ge 0\}$ be moved by random
stirring on $\Zz$, with $Y(0)=0$ and $Y'(0)=1$ and let $p(u)$ be
the probability that at time $u$ the two particles are exchanged,
i.e., $p(u)=\Pr[Y(u) >Y'(u)]$, and note that this is the same as
requiring the particles to have been exchanged an odd number of
times. For all $t>0$, let
\begin{equation}
\label{def-f}
 f(t):= \sum_{x=1}^{\infty}\int_0^t
\Pr[T^x \in ds]p(t-s)
\end{equation}

\begin{theorem} Let $t > 0$. Then $n^{-1} D^n_{nt} \to_p f(t)$
as $n \to \infty$ where $f$ is the function defined by
(\ref{def-f}). $f(t)$ is infinitely differentiable, and moreover
it has the asymptotic behavior
$$
\lim_{t\to \infty} \frac{f(t)}{\sqrt{t}}
= \half E\left( \max_{0\le s\le 1} B_{4s} \right) = \sqrt{\frac{2}{\pi}}
$$
where $B_t$ is a standard Brownian motion.
\label{adjsmt}
\end{theorem}

\noindent
To check the constant recall that
$$
P\left( \max_{0\le s\le 1} B_{4s} \right) = 2 P( B_4 > x )
$$
so integrating gives
\beqax
\half E\left( \max_{0\le s\le 1} B_{4s} \right) & = & \int_0^\infty P( B_4 > x ) \, dx
= EB_4^+ = 2 EB_1^+ \cr
& = & \frac{2}{\sqrt{2\pi}} \int_0^\infty x e^{-x^2} \, dx = \sqrt{\frac{2}{\pi}}
\eeqax

The next result looks at the distance of the random walk at times of order $n^3$,
i.e., when each particle has moved of order $n^2$ times, and hence has
a significant probability of hitting a boundary. Let $p_t(u,v)$ denotes the transition
function of $\bar B$, a one-dimensional Brownian motion
run at speed 2 reflecting at 0 and 1.

\begin{theorem} Let $t>0$.
$$
\frac1{n^2} D^n_{n^3t} \to_p
\int_0^1du\int_u^1dv\int_0^1p_t(u,x)dx\int_0^y
p_t(v,y)dy=\Pr[ \bar B_1(t)> \bar B_2(t)]
$$
where $\bar B_1$ and $\bar B_2$ are independent copies of $\bar B$
started uniformly on $0\le \bar B_1(0) < \bar B_2(0)\le 1$
evolving independently. \label{adjlgt}
\end{theorem}

In between the two extremes we have a simple behavior.

\begin{theorem}
\label{adjint}
Let $s=s(n)$ with $s \to \infty$ and $s/n^2 \to 0$. Then
$$
\frac{1}{n\sqrt{s}} D^n_{ns} \to_p \sqrt{\frac{2}{\pi}}
$$
\end{theorem}

\mn Recently, Angel et al. \cite{sorting} have also used the
simple exclusion process to analyze a process on the Cayley graph
of the symmetric group generated by adjacent transpositions, but
this time in the context of sorting networks.

\subsection{Random walk on a random 3-regular graph}

A \emph{3-regular} graph is a graph where all vertices have degree
equal to 3. To construct a random 3-regular graph we suppose $n$
is even, and use the approach of Bollob\`as and de la Vega
\cite{bb-vega82} (see also Bollob\`as \cite{bb88}). Expand each
vertex $i$ into 3 ``mini-vertices" $3i,\ 3i+1$ and $3i+2$, and
consider a random matching $\sigma(j)$ of the $3n$ mini-vertices.
A random 3-regular graph $G_n$ is then obtained by collapsing back
the $n$ groups of 3 mini-vertices into $n$ vertices while keeping
the edges from the random matching. We may end up with self-loops
or multiple edges, but with a probability that is positive
asymptotically, we do not, so the reader who wants a neat graph
can condition on the absence of self-loops and multi-edges.

Departing from our choices in the previous example, we consider
the discrete time random walk $\hat X^n_k$, $k\ge 0$, that jumps
from $j$ to $[\sigma(3j+i)/3]$ where $i$ is chosen at random from
$\{0,1,2\}$. (We have used this definition since it works if there
are self-loops or multiple edges.) Let $\hat D^n_t$ be the
distance from the starting point at time $t$.

\begin{theorem}
\label{rgt}
For fixed $t>0$
$$
\frac{\hat D^n_{[t\log_2 n]}}{\log_2 n} \to_p \min\left( \frac{t}{3},1\right)
$$
\end{theorem}

An intuitive description of a random 3-regular graph, as seen from
vertex 1, can be given as follows. Grow the graph by successively
adding vertices adjacent to the current set. Branching process
estimates will show that as long as the number of vertices
investigated is $O(n^{1-\ep})$, this portion of the graph looks
very much like a regular tree in which each vertex has 2 edges
going away from the root and 1 leading back towards the root.
Thus, until the distance of $\hat X_n$ from $\hat X_0$ is $\ge
(1-\epsilon)\log_2 n$, $\hat D^n_k$ evolves like a (2/3,1/3)
biased random walk on the nonnegative integers, with transition
probabilities $p(x,x+1)=2/3$ and $p(x,x-1)=1/3$, and reflection at
0. After $k$ moves we expect this walk to be at distance $k/3$. On
the other hand, once the walk reaches a distance corresponding to
the diameter of the graph, $\log_2 n$ by Bollob\`as and de la Vega
\cite{bb-vega82}, or Theorem 2.13 in Worwald \cite{wormwald}, it
should remain at this level. Indeed, it cannot go any further,
since this is the diameter. On the other hand the tree structure
below makes it hard for it to come down back toward the root.

\mn {\bf Open Problem.} The techniques developed for the random
walk on a 3-regular graph should be useful when dealing with
random walk on the giant cluster of a Erd\H{o}s-R\'enyi random
graph with $p=c/n$ and $c>1$, which locally has the geometry of a
``Poisson mean $c$ Galton-Watson tree". We conjecture that the
random walk exhibits a phase transition like the one in Theorem
\ref{rgt} but with a different constants in place of 3 and 1 on
the right-hand side. One technical problem is that the diameter is
strictly larger than the average distance between points $\log
n/(\log c)$, see Chung and Lu \cite{chung-lu}, so we don't have
the easy upper bound.

\subsection{Other random walks}

We view the analysis of the examples above as the start of a more
complete investigation of what are the possible behaviors for the
distance of a random walk. There are many other interesting
examples to consider. For instance, Fulman \cite{fulman} has
studied the evolution of the distance for the
Gilbert-Shannon-Reeds riffle shuffle. For this shuffling method
(which, we emphasize, is \emph{nonreversible}), Bayer and Diaconis
\cite{bayer-diaconis} proved an explicit formula for the
distribution of $\sigma$ after $r$ shuffles. In particular from
this formula it follows that the riffle-shuffle distance of a
permutation is given by
$$
d(\sigma)=\lceil\log_2 (\text{Des}(\sigma)+1)\rceil \text{ where
Des($\sigma$)}= \#\{1\le i\le n-1: \sigma(i)>\sigma(i+1)\}
$$
Des($\sigma$) is called the number of descents of $\sigma$. The
main result of Fulman \cite{fulman} is
\begin{theorem}
After $r= \log_2(\alpha n)$ shuffles
$$
\frac1n E(\text{Des}(\sigma_{r}))\to \alpha-\frac1{e^{1/\alpha}-1}
$$
at least if $\alpha>1/(2\pi)$.
\end{theorem}
This says that for this range of $r$ the walk is already in a
sublinear regime. In particular, as $\alpha \to \infty$ we get
that $d(\sigma_r)\sim \log_2 n$ since the expression
$$
\alpha-\frac1{e^{1/\alpha}-1} = \alpha - \alpha \cdot \frac{ 1 } { 1 + (1/\alpha)/2 + \cdots } \to 1/2
$$
This is not surprising since this is the diameter of the graph. It
is not clear at this point whether this formula also holds for
smaller values of $\alpha$, although it is tempting to let $\alpha
\to 0$ and get that for small values of $\alpha$ the walk is
``almost" linear (the fraction term with the exponential is much
smaller than the other term).

\section{Random adjacent transpositions}

Let $X_t$, which we write for now on without the superscript $n$,
to be the continuous time walk on permutations of $\{1, 2, \ldots n\}$
in which at rate 1 we pick a random $1\le i \le n-1$ and exchange the
values of $X_t(i)$ and $X_t(i+1)$. As indicated in (\ref{inv})
the distance from a permutation $\sigma$ to the identity is
$\adj(\sigma)=\#\{i<j:\ \sigma(i)>\sigma(j)\}$,
the number of inversions of $\sigma$.

\subsection{Small times}

The reflecting boundaries at 1 and $n$ are annoying complications, so the first thing we will do
is get rid of them. To do this and to prepare for the variance estimate
we will show that if $i<j$ are far apart then the probability
$X_t(i) > X_t(j)$ is small enough to be ignored. Let $P^{[a,b]}$ be
the probabilities for the stirring process with reflection at $a$ and $b$,
with no superscript meaning no reflection.

\begin{lemma}
$ P^{[1,n]} (X_{nt}(i) > X_{nt}(j))
\le 8 P( X_{nt}(0) > (j-i)/2 )$
\label{refcomp}
\end{lemma}

\begin{proof} A simple coupling shows
\beqax
P^{[1,n]} (X_{nt}(i) > X_{nt}(j))
& \le &P^{[0,j-i]} (X_{s}(0) > X_{s}(j-i) \hbox{ for some $s\le nt$}) \cr
& \le & 2 P^{[0,\infty)}\left( \max_{0 \le s \le nt} X_s(0) > (j-i)/2 \right)
\eeqax
Using symmetry and then the reflection principle, the last quantity is
$$
\le  4 P\left( \max_{0 \le s \le nt} X_s(0) > (j-i)/2 \right)
\le  8 P( X_{nt}(0) > (j-i)/2 )
$$
which completes the proof.
\end{proof}

Since the random walk on time scale $nt$ moves at rate 2,
$$
E \exp(\theta X_{nt}(0))  =  \sum_{k=0}^\infty e^{-2t} \frac{(2t)^k}{k!}
\left( \frac{ e^\theta + e^{-\theta} }{2} \right)^k
 =  \exp( - 2t + t (e^\theta + e^{-\theta}) )
$$
Using Chebyshev's inequality, if $\theta > 0$ \beq P( X_{nt}(0) >
x ) \le \exp( - \theta x + t [ e^\theta + e^{-\theta} - 2 ])
\label{ldbd} \eeq Taking $\theta=1$
$$
P( X_{nt}(0) > x ) \le C_t e^{-x} \quad\hbox{where}\quad C_t = \exp( (e + e^{-1} - 2 )t )
$$

When $x=3 \log n$ the right-hand side is $C_t n^{-3}$, so using Lemma
\ref{refcomp}, for fixed $t$ it suffices to consider ``close pairs" with
$0 < j-i < 6 \log n$. The number of close pairs with $i \le 3 \log n$
or $j > n - 3 \log n$ is $\le 36 \log^2 n$, so we can ignore these
as well, and the large deviations result implies that it is enough to consider
random stirring on $\Zz$.

We are now ready to prove the first conclusion in Theorem
\ref{adjsmt}: if $t > 0$ then as $n \to \infty$ \beq \frac1n
D^n_{nt} \to f(t)= \sum_{x=1}^{\infty}\int_0^t \Pr[T^x \in
ds]p(t-s) \label{th2a} \eeq

\begin{proof}[Proof of $(\ref{th2a})$]
It is clear from the Markov property that if $X_t(i)$ and $X_t(j)$
are moved by stirring on $\Zz$ then
$$
P( X_t(i) > X_t(j) ) = \int_0^t \Pr[T^{j-i} \in ds]p(t-s)
$$
With the large deviations bound in (\ref{ldbd}) giving us domination we can
pass to the limit to conclude
$$
\frac{1}{n} \sum_{1\le i < j \le n } P( X_t(i) > X_t(j) )
\to \sum_{x=1}^\infty \int_0^t \Pr[T^{x} \in ds]p(t-s)
$$

To prove convergence in probability let
$$
\xi_{i,j} =  1_{( X_t(i) > X_t(j) )} - P( X_t(i) > X_t(j) )
$$
By remarks above it suffices to consider the sum over
$1\le i < j \le n$ with $0< j-i< 6\log n$,
$i \ge 3 \log n$ and $j \le n - 3 \log n$ which we denote
by $\Sigma^*$, and if $i' > j + 6\log n$ then
$$
E\xi_{i,j} \xi_{i'j'} \le 4C_t n^{-3}
$$
since the random variables have $|\xi|\le 1$ and will be
independent unless some random walk moves by more than $3\log n$
in the wrong direction. From this it follows that
$$
E \left( \Sigma^* \xi_{i,j} \right)^2 \le n \cdot (6\log n)^3 + 4C_t n^{-3} (n \cdot 6 \log n)^2
$$
and the result follows from Chebyshev's inequality.
\end{proof}

The remaining detail is to show that $f$ is smooth and that as
$t\to \infty$, \beq \lim_{t \to \infty} f(t)/\sqrt{t} = \half \Ee
\left( \max_{0\le s\le 1} B_{4s} \right) \label{th2b} \eeq where
$B_{\cdot}$ is a standard Brownian Motion. The fact that $f$ is
infinitely differentiable follows easily from repeated use of
Lebesgue's theorem and the fact that both $p(u)$ and $dP(T^x\in
du)/du$ are infinitely differentiable smooth functions. This is
itself easily checked: for instance, if $q_j$ is the probability
that a simple random walk in discrete time started at 0 hits $x$
in $j$ steps, then $T^x=\sum_{j=x}^{\infty} q_j
\text{Gamma}(j,4)$, so $T^x$ has a smooth density. A similar
argument also applies for the function $p(u)$.

\begin{proof}[Proof of $(\ref{th2b})$] The result follows easily from two simple lemmas.

\begin{lemma}
$p(t) \to 1/2$ as $t \to \infty$.
\end{lemma}

\begin{proof} Each time there is jump when the particles $Y$ and $Y'$ are adjacent, they have a
probability $1/3$ of being exchanged the next step. So,
conditionally on the number of such jumps $N$, the number of
actual swaps between $Y$ and $Y'$ is Binomial$(N, 1/3)$. Now, $Y>Y'$ if
and only if the number of times they are swapped is odd. Hence the
lemma follows from the two observations : (i) As $t \to \infty$,
the number of jumps while they are adjacent to each other $\to \infty$,
and (ii) as $N \to \infty $, $\Pr[\hbox{Binomial}(N,p)\text{ is odd}] \to 1/2$
for any given $0<p<1$. For (i), observe that the discrete-time
chain derived from $\{|Y_t-Y'_t|-1, t \ge 0\}$ is a reflecting
random walk on $\{0, 1, \ldots\}$, and therefore visits 0
infinitely many times. (ii) is an easy fact for Bernoulli random
variables.
\end{proof}

\begin{lemma}
$$
\frac1{\sqrt{t}}\sum_{x=1}^{\infty} \Pr[T^x \in (t-\log t, t)] \to 0
$$
\label{Txsumbd}
\end{lemma}

\begin{proof}
The random walk can only hit a new point when it jumps so
\beqax
t^{-1/2}E\left( \sum_{x=1}^{\infty}\indic{T^x\in (t-\log t,t)} \right)&\le&
t^{-1/2}E(\#\text{jumps of the random walk in} (t-\log t,t))\\
&\le & t^{-1/2}\cdot (4\log t) \to 0
\eeqax
since jumps occur at rate 4.
\end{proof}

It is now straightforward to complete the proof. Let $\eps >0$.
Fix $T$ large enough so that $|p(t)-1/2|\le \eps$ as soon as $t
\ge T$. Then by Lemma \ref{Txsumbd}, for $t \ge T':=e^T$, letting
$W_t$ be a simple random walk on $\Zz$ in continuous time jumping
at rate 1, \beqax
t^{-1/2}f(t)&=&t^{-1/2}\sum_{x=1}^{\infty}\int_0^{t-\log
t}\Pr[T^x\in ds]p(t-s) +o(1)
\\
&\le& \left(\half + \eps\right)t^{-1/2}\sum_{x=1}^{\infty}\Pr[T^x <t-\log t] +o(1)
\\
&\le &\left( \half + \eps
\right)t^{-1/2}\sum_{x=1}^{\infty}\Pr\left(\max_{s\le t- \log t}
W_{4s} >x \right)
\\
&\to& \left(\half+\eps\right) \Ee\ \max_{s\le 1} B_{4s}. \eeqax
by Donsker's theorem.  The other direction
$\liminf_{t\to\infty} t^{-1/2}f(t) \ge (1/2 - \eps)
\Ee \max_{s\le 1} B_{4s}$ can be proved in the same way.
\end{proof}

\subsection{Large times}

Our next goal is to prove that if $t>0$ then
\beq
\frac1{n^2} D^n_{n^3t} \to
\int_0^1du\int_u^1dv\int_0^1p_t(u,x)dx\int_0^y
p_t(v,y)dy=\Pr[\bar B_1(t)> \bar B_2(t)]
\eeq
in probability and where $\bar B_1$ and $\bar B_2$ are two reflecting
Brownian motions run at speed 2 started uniformly on $0\le \bar B_1(0) < \bar B_2(0)\le 1$ and evolving independently.

\begin{proof} We first show that the expected value converges.
The first step is to observe that the rescaled random walks
$X_{n^3t}(i)/n$, $t\ge 0$ converge to reflecting Brownian Motion
on $[0,1]$. Indeed, Durrett and Neuhauser \cite[(2.8)]{dn94}
showed that for fixed $i<j$, the rescaled pair of random walks
converge to two independent Brownian Motions. They did this on
$\Zz$ but the proof extends in a straightforward way to the
current setting. Their proof shows that if $i/n \to x$ and $j/n
\to y$ we have
$$
\Pr[X_{n^3t}(i)>X_{n^3t}(j)] \to \Pr_{x,y}[\bar B_1(t)> \bar B_2(t)]
$$
This implies that the convergence occurs uniformly on the compact set so
$$
\frac1{n^2}
\Ee D^n_{n^3t} = \frac1{n^2}\sum_{i<j}\Pr[X_{n^3t}(i)>X_{n^3t}(j)]
\to\Pr[\bar B_1(t)>\bar B_2(t)]
$$

To get the convergence in probability, we use second moment estimates.
Let $A_{i,j} = \{ X_{n^3t}(i)>X_{n^3t}(j) \}$.
$$
\Ee\left(\frac1{n^2} D^n_{n^3t}\right)^2 =
\frac1{n^4}\sum_{i<j}\sum_{k<l}\Pr[A_{i,j} \cap A_{k,\ell}]
$$
The first step is to observe that there are only $O(n^3)$ terms in
which two of the indices are equal so these can be ignored. When
the four indices are distinct we can again apply Durrett and
Neuhauser's \cite{dn94} results to the 4-tuple of random walks
$(X(i),X(j),X(k),X(l))$, to conclude that if $i/n \to x$, $j/n \to
y$, $k/n \to x'$ and $l/n \to y'$
$$
\Pr[A_{i,j} \cap A_{k,\ell}] \to \Pr_{x,y}[B_1(t)>B_2(t)]
\Pr_{x',y'}[B_1(t)>B_2(t)]
$$
From this it follows that
$$
\Ee\left(\frac1{n^2} D^n_{n^3t}\right)^2 - \left( \Ee\frac1{n^2} D^n_{n^3t}\right)^2  \to 0
$$
In other words, the variance of
$n^{-2}D^n_{n^3t}$ is asymptotically 0, and applying
Chebyshev's inequality, we get the convergence in probability
to the limit of the means.
\end{proof}

\subsection{Intermediate regime}

The proof of Theorem \ref{adjint} is a hybrid of the two previous proofs.
We first truncate to show that it suffices to consider $i<j$ close together
and far from the ends, then we compute second moments. We begin with a
large deviations result:

\begin{lemma}
For all $x>0$ and $t>0$ then
$$
P( X_{nt}(0) > x ) \le \exp( - x^2/8et ) + \exp(-x\ln(2)-2t)
$$
\label{ld2}
\end{lemma}

\begin{proof} First assume $x\le 4e t$. From (\ref{ldbd}) we have
$P( X_{nt}(0) > x ) \le \exp( - \theta x + t [ e^\theta +
e^{-\theta} - 2 ] )$. When $0 < \theta < 1$ \beqax e^\theta +
e^{-\theta} - 2 & = & 2 \left[ \frac{\theta^2}{2} +
\frac{\theta^4}{4!} + \frac{\theta^6}{6!} + \cdots  \right] \cr &
\le & \theta^2 \left[ 1 + \frac{\theta^2}{2^2} +
\frac{\theta^4}{2^4} + \cdots  \right] = \frac{\theta^2}{1-
\theta^2/4} \le \frac{4\theta^2}{3} \eeqax and by continuity this
is valid also when $\theta=1$. Taking $\theta = x/4et$ which is
$\le 1$ by assumption
\begin{equation}
\label{ld21} P( X_{nt}(0) > x ) \le \exp\left(  - \frac{x^2}{4et}
+ t \cdot \frac{4}{3} \left(\frac{x}{4et} \right)^2 \right) \le
\exp(-x^2/8et)
\end{equation}
When $x>4et$, remark that $P(X_{nt}(0)>x)$ is smaller than the
probability that a Poisson random variable with mean $2t$ is
greater than $x$. Thus for any $\theta>0$ this is by Markov's
inequality smaller than $\exp(-\theta x + 2t(e^\theta-1))$. This
is optimal when $e^\theta= x/2t$, in which case we find that
\begin{equation}\label{ld22}
P(X_{nt}(0)>x)\le \exp(-x\ln(x/2t)+x-2t)\le \exp(-x\ln(2)-2t)
\end{equation}
since $x\ge 4et$. Equations (\ref{ld21}) and (\ref{ld22}) give us
two bounds valid in different regions, so by summing them we get a
bound that is everywhere valid, and this concludes the proof.
\end{proof}

\begin{proof}[Proof of Theorem \ref{adjint}]
By assumption we can pick $K_n \to\infty$ so that $K_n^2
\sqrt{s}/n \to 0$. By Lemma \ref{refcomp},
\begin{eqnarray*}
\frac{1}{n\sqrt{s}} \sum_{i,j > i + K_n\sqrt{s}} P^{[1,n]}(
X_{ns}(i)
> X_{ns}(j) ) &\le & 8 \frac1{\sqrt{s}} \sum_{x=K_n\sqrt{s}}^{\infty}
P(X_{ns}(0)>x/2)\\
&= & 8 \int_{K_n}^{\infty}P(X_{ns}(0)>\lfloor x\sqrt{s}/2\rfloor )
dx
\end{eqnarray*}
Applying Lemma \ref{ld2} it follows that
$$
\frac{1}{n\sqrt{s}} \sum_{i,j > i + K_n\sqrt{s}} P^{[1,n]}(
X_{ns}(i)
> X_{ns}(j) ) \to 0
$$
Letting $I_{i,j}$ be the indicator of $\{ X_{ns}(i)
> X_{ns}(j) \}$ it follows that
$$
\frac{1}{n\sqrt{s}} \sum_{i,j > i + K_n\sqrt{s}} I_{i,j} \to 0 \quad\hbox{in probability}
$$
i.e., we can restrict our attention to close pairs. Once we do this, we can eliminate
ones near the ends since
$$
\frac{1}{n\sqrt{s}} \sum_{i< K_n \sqrt{s},j \le i + K_n\sqrt{s}} I_{i,j}
\le \frac{ (K_n \sqrt{s})^2 }{ n \sqrt{s} } \to 0
$$
by assumption. In a similar way we can eliminate $j > n - K_n \sqrt{s}$.

It follows that it is enough to consider random stirring on $\Zz$.
The result of Durrett and Neuhauser \cite{dn94} implies that if $s
\to \infty$, $i \ge K_n \sqrt{s}$, $j \le n - K_n \sqrt{s}$ and
$(j-i)/\sqrt{s} \to x$ then
$$
EI_{i,j} \to \frac{1}{2} P\left( \max_{0\le t\le 1} B_{4t} > x \right)
$$
where the right-hand side is 0 if $x=\infty$. Writing $\Sigma^*$
again for the $i,j$ with $i \ge K_n \sqrt{s}$, $j \le n - K_n
\sqrt{s}$ and $0 < j-i < K_n \sqrt{s}$, and using the domination
that comes from Lemma \ref{ld2} it follows that
$$
\frac{1}{n\sqrt{s}} \Sigma^* EI_{i,j} \to \frac{1}{2} E \max_{0 \le t \le 1} B_{4t}
$$

The next step is to compute the second moment. The number of terms
with one index in $i<j$ equal to one of $k<l$ with both pairs
close is $\le n (K_n \sqrt{s})^2$, which when divided by
$(n\sqrt{s})^2$ tends to 0. The result of Durrett and Neuhauser
\cite{dn94} implies that terms in which all four indices are
different are asymptotically uncorrelated. We remark that in Lemma
\ref{ld2} we can also get an upper-bound on $P(X_{nt}(0)>x)^{1/2}$
by summing the square-roots of the two terms in (\ref{ld21}) and
(\ref{ld22}) since only one of them applies in a given region.
This and Cauchy-Schwartz's inequality provide the justification
for the passage to the limit:
$$
\frac{1}{(n\sqrt{s})^2} \Sigma^*_{i<j} \Sigma^*_{k<\ell}
E(I_{i,j}I_{k,\ell}) \to \left( \frac{1}{2} E \max_{0 \le t \le 1}
B_{4t} \right)^2
$$
and the rest of the argument is the same as in Theorem
\ref{adjlgt}.
\end{proof}

\section{Random walk on a random 3-regular graph}

Let $G_n$ be a random 3-regular graph constructed as in the introduction
and let $X_k$ be the discrete time random walk on $G_n$, where
for simplicity we drop both the superscript $n$ and the hat
to indicate discrete time. We assume that $X_0=1$
and write $D^n_k$ for the graph distance from $X_k$ to $X_0$.
Our goal is to prove Theorem 5, that is, for fixed $t>0$
\beq
\frac{D^n_{[t\log_2n]}}{\log_2 n} \to_p \min\left( \frac{t}{3},1\right)
\label{th5}
\eeq

\subsection{Proof for the subcritical regime}

Let $v$ be a vertex at distance $l$ from the root. We say that $v$
is a ``good" vertex if it has two edges leading away from the root
(at distance $l+1$) and one leading back to distance $l-1$.
Otherwise we say that $v$ is a ``bad" vertex. Let $B(l)$ be the
set of all bad vertices at distance $l$.

\begin{lemma}
\label{badv} Let $2\le v\le n$ be a vertex distinct from the
root. Given that $v$ is at distance $l$ from the root,
$P(v\in B(l))\le 2i/n$ where $i=2^l$.
\end{lemma}

\begin{proof} First consider the event
that $v$ has an edge leading to some other vertex at distance $l$.
Since it is at distance $l$, it must have at least one edge
leading backwards, so there are only two other edges left. In
particular there are at most $2^l=i$ vertices at distance $l$. In
$G_n$ those $i$ vertices at distance $l$ correspond to $2i$
unpaired mini-vertices, so the probability of a connection
sideways to another vertex at distance $\ell$ is smaller than $2i/3n$.

When $v$ has two edges leading forward, the probability that
one of its children is connected to another vertex from level $l$
is also smaller than $2i/3n$ since there are at most $2i$ edges
leading to level $l+1$. Since $v$ has at most 2 children, this
gives a probability of at most $4i/3n$. Combining this with the
estimate above gives $2i/3n +4i/3n=2i/n$.
\end{proof}

A simple heuristic now allows us to understand that with high probability the random
walk will not encounter any vertices as long as we are in the subcritical regime.
Before we encounter a bad vertex, the distance is a (2/3,1/3) biased random walk and hence spends
an average of 2 steps at any level. Hence, the
expected number of bad vertices encountered until time distance
$(1-\eps)\log_2 n$ is smaller than
$$
\sum_{l=1}^{(1-\eps)\log_2 n} 2 \frac{2^l}n=O(n^{-\eps})\to 0
$$
To prove this rigorously, let $A_k$ denote the event that by time $k$ the
random walk has never stepped on a bad vertex up to time $k$.

\begin{lemma}
As $n\to\infty$, $P( A_{3(1-\eps)\log_2 n} ) \to 1$.
\label{Akbd}
\end{lemma}

\noindent
On this event, for each $1\le j\le 3(1-\eps)\log_2 n$, $X_j$ has probability 2/3 to move away
from the root and 1/3 to move back towards the root,
and the first part of Theorem \ref{rgt} follows easily.

\begin{proof} By Lemma \ref{badv} the probability that some vertex within
distance $L$ of 1 is bad is
$$
\le \sum_{\ell=1}^L 2^\ell \frac{2 \cdot 2^\ell}{n}
\le \frac{2}{n} \cdot \frac{2^{2L}}{1 - 1/4} \to 0
$$
if $L = (1/3) \log_2 n$.

Since for each vertex there are at most two edges leading out
and one leading back, the distance from the starting point is bounded
above by a (2/3,1/3) biased random walk. Standard large deviations arguments
imply that there are constants $C$ and $\alpha$ depending on $\rho$ so that
\begin{equation}
\label{upbeld}
P(d(X_k)>\rho k)\le Ce^{-\alpha k}
\end{equation}
Summing from $k=L$ to $\infty$, we see that with high probability
$d(A_k) \le \rho k$ for all $k \ge L$.

When this good event occurs for $k \ge L$, it follows from Lemma
\ref{badv} that
$$
P(A_{k+1}) \ge P(A_k)\left(1- 2\frac{2^{\rho k}}n\right)
           \ge \prod_{j=L}^k \left(1-2\frac{2^{\rho j}}n\right)
$$
Taking the logarithm, we have for large $n$
\beqax
\log P(A_{k+1}) &\ge& \sum_{j=L}^{k}\log\left(1- 2\frac{2^j}{n}\right)\\
                &\ge& -4\sum_{j=1}^{k}\frac{2^j}{n}
            \ge -\frac{4}{1 - 2^{-\rho}} \cdot \frac{2^{k\rho}}{n}
\eeqax
We want to take $k=3(1-\eps)\log_2 n$. By choosing $\rho$ close enough to $1/3$ so that
$3\rho(1-\eps)<1$, we have $2^{k\rho}/n = n^{-\alpha}$ with $\alpha>0$ which proves the desired result.
\end{proof}

\subsection{Proof for the supercritical regime}

Here we wish to prove that if $k=t\log_2 n$, with $t>3(1-\eps)$, then
$d(X_k)\approx \log_2 n$. As already noted, this is the diameter
of $G_n$ so all we have to prove is that once it reaches this
distance it stays there. To do this we let
$$
L(a,b):=\{2\le v\le n:d(v)\in[a\log_2n,b\log_2n]\}
$$
and consider $L(1-\eps,1-\delta)$.

Intuitively, this strip consists of about $n^{1-\eps}$ trees, each
with at most $n^{\eps-\delta}$ vertices. However, there are
sideways connections between these trees so we have to be careful
in making definitions. Let $v_1,\ldots,v_m$ be the $m$ vertices at
level $(1-\eps)\log_2n$. For $j=1,\ldots,m$ if $v\in
L(1-\eps,1-\delta)$, we say that $v\in T_j$ if $v_j$ is the
closest vertex to $v$ among $v_1,\ldots,v_m$.

To estimate the number of sideways connections (i.e., edges
between vertices $v$ and $v'$ in different $T_j$'s), we use:

\begin{lemma}
\label{swcon} The number of subtrees that $T_j$ is
connected to is dominated by a branching process with offspring
distribution Binomial$(n^{\eps-\delta},n^{-\delta})$.
\end{lemma}

\begin{proof}
Each tree to which we connect requires a bad connection (i.e., one
of the two possible errors in Lemma \ref{badv}). Suppose we
generate the connections sequentially. The upper bound in Lemma
\ref{badv} holds regardless of what happened earlier in the
process, so we get an upper-bound by declaring each vertex at
level $l$ bad independently with probability $2i/n$ with $i=2^l$,
so this probability is at most $ n^{-\delta}$. Since there are at
most $n^{\eps-\delta}$ vertices in a given subtree, the lemma
follows immediately.
\end{proof}

\begin{lemma}
\label{badcl} If $\delta>\eps/2$ then there exists some
$K=K(\eps,\delta)>0$ such that,
$$
P(\hbox{there is a cluster of trees $T_j$ with more than $K$ bad vertices})\to 0
$$
\end{lemma}
\begin{proof}
The worst case occurs when each bad connection in a tree leads to a new one. Let
$$
X\overset{d}{=}\text{Bin}(n^{\eps-\delta},n^{-\delta})
$$
be the offspring distribution of the branching process of the
previous Lemma. In particular
$$
E(X)=O(n^{\eps-2\delta})\to 0
$$
Let $c=n^{\eps-2\delta}$, and let $N=n^{\eps-\delta}$ be the
total number of vertices in $T_j$, so
$X\overset{d}{=}\text{Binomial}(N,c/N)$.

Lemma \ref{swcon} follows from a simple evaluation of the tail of the total progeny
$Z$ of a branching process with offspring distributed as $X$. To
do this, we let
\beqax
\phi_N(\theta) & = & e^{-\theta} \sum_{k=0}^{N} \binom{N}{k}
\left( \frac{c}{N} \right)^k
\left( 1 - \frac{c}{N} \right)^{N-k} e^{\theta (k-1)} \\
& = & e^{-\theta} \left( 1- \frac{c}{N} + \frac{c}{N} e^\theta \right)^{N}
\eeqax
be the moment generating function of $X- 1$.
Let $S_k$ be a random walk that takes steps with this distribution and $S_0=1$. Then
$\tau = \inf\{ k: S_k = 0\}$ has the same distribution as $Z$. Let
$R_k = \exp(\theta S_k)/\phi_N(\theta)^k$. $R_k$ is a nonnegative
martingale. Stopping at time $\tau $ we have $e^{\theta} \ge
E(\phi_N(\theta)^{-\tau})$. If $\phi_N(\theta)<1$ it follows that
$$
P(\tau \ge y ) \phi_N(\theta)^{-y} \le E[\phi_N(\theta)^{-\tau}]
\le e^{\theta}
$$
Using $\phi_N(\theta) \le e^{-\theta} \exp(c(e^\theta - 1))$ now
we have
$$
P(\tau \ge y ) \le e^{\theta} \left(e^{-\theta}\exp(c(e^\theta -
1))\right)^y
$$
To optimize the bound we want to minimize $c(e^\theta - 1) -
\theta$. Differentiating this means that we want $ce^\theta - 1 =
0$ or $\theta = -\log(c)$.  Plugging this and recalling that
$\tau$ and $Z$ have the same distribution we have
$$
P(Z \ge y ) \le \frac1c \exp( - (c-1 - \ln c) y )
$$
Substituting $c=n^{-\alpha}$ with $\alpha =2\delta-\eps$, we find
that
$$
P(Z\ge y)\le n^{\alpha}\exp(y(1-\alpha \log(n)))
$$
Since there are $m\le n^{1-\eps}$ trees to start with, the
probability that one of them has more than $y$ trees in its
cluster is smaller than
$$
n^{1-\eps} n^{\alpha}\exp(y(1-\alpha \log n))
$$
so if
$$
y>\frac{\alpha+1-\eps}{\alpha}:=K(\delta,\eps)
$$
then the probability than one cluster contains more than $y$ trees
tends to 0. This implies that with probability 1 asymptotically,
no cluster of trees has more than $K$ bad vertices, since the
branching process upper-bound is obtained by counting every bad
vertex as a sideways connection.
\end{proof}

With Lemma \ref{badcl} established the rest is routine. In each
cluster of trees there is a stretch of vertices of length $\ge a
\log_2n$ where $a=(\ep-\delta)/(K+1)$ with no bad vertices. The
probability of a downcrossing of such a strip by a (2/3,1/3)
random walk is $\le (1/2)^{a \log_2 n} = n^{-a}$ so the
probability of one occurring in $n^{a/2}$ time steps tends to 0.

%\section*{References}


\begin{thebibliography}{99}

\bibitem{sorting} O. Angel, A. Holroyd and D. Romik. Random sorting networks. In preparation.

\bibitem{bayer-diaconis} D. Bayer and P. Diaconis (1992). Trailing the dovetail
shuffle to its lair. \emph{Ann. Probab.}, 2, 294-313.

\bibitem{bd} N. Berestycki and R. Durrett, (2006). A phase
transition in the random transposition random walk. \emph{Probab.
Theory Rel. Fields}, to appear.

\bibitem{bb-book} B. Bollob\'as (1985). \emph{Random graphs}.
Academic Press, London.

\bibitem{bb88} B. Bollob\'as (1988). The isoperimetric number of a
random graph, \emph{European Journal of Combinatorics}, 9,
241-244.

\bibitem{bb-vega82} B. Bollob\'as and F. de la Vega (1982). The
diameter of random regular graphs. \emph{Combinatorica}, 2,
125-134

\bibitem{chung-lu} F.K. Chung and L. Lu (2001). The diameter of sparse random
graphs. \emph{Adv. Appl. Math.} 26, 257-279.


\bibitem{diaconis} P. Diaconis (1988). \emph{Group representation in
Probability and Statistics}, Institute of Mathematical Statistics
Lecture Notes, Vol. 11.

\bibitem{dn94} R. Durrett and C. Neuhauser (1994). Particle systems
and reaction-diffusion equations. \emph{Ann. Prob.}, Vol. 22, No.
1, 289-333.

\bibitem{eriksen} N. Eriksen (2005). Expected number of inversions
after a sequence of random adjacent transpositions - an exact
expression. \emph{Discrete Mathematics}

\bibitem{erikkson+} H. Eriksson, K. Erikkson, and J. Sj\"ostrand
(2000). Expected number of inversions after $k$ random adjacent
transpositions. In D. Krob, A.A. Mikhalev, A.V. Mikhalev, eds.
\emph{Proceedings of Formal Power Series and Algebraic
Combinatorics}, Springer-Verlag (2000) 677-685

\bibitem{fulman} J. Fulman (2005). Stein's method and minimum
parsimony distance after shuffles. \emph{Electr. J. Probab.} 10,
901--924.

\bibitem{wormwald} N.C. Wormald (2005). Models of random regular
graphs (survey). Available at \linebreak {\tt
http://www.ms.unimelb.edu.au/$\sim$nick/papers/regsurvey.pdf}

\end{thebibliography}
\end{document}